\documentclass[12pt]{article}
\usepackage{amsmath,color}
\usepackage{graphicx}
\usepackage{enumerate}
\usepackage{natbib}
\usepackage{url} 
\usepackage{multirow}
\usepackage{multicol}
\usepackage{nicefrac}
\usepackage{subfigure}
\usepackage{algorithm,algpseudocode}
\newcommand{\blind}{1}

\newcommand{\bmu}{\mbox{\boldmath $\mu$}}
\newcommand{\bphi}{\mbox{\boldmath $\phi$}}

\newcommand{\bSigma}{\mbox{\boldmath $\Sigma$}}
\newcommand{\bepsilon}{\mbox{\boldmath $\epsilon$}}

\newcommand{\bbeta}{\mbox{\boldmath $\beta$}}

\newcommand{\bGamma}{\mbox{\boldmath $\Gamma$}}

\newcommand{\yp}{\mathbf{y}}
\newcommand{\xp}{\mathbf{x}}
\newcommand{\y}{\mathbf{y}}
\newcommand{\Y}{\mathbf{Y}}

\newcommand{\zp}{\mathbf{z}}
\newcommand{\Z}{\mathbf{Z}}
\newcommand{\R}{\mathbf{R}}
\newcommand{\X}{\mathbf{X}}

\newcommand{\sa}{\mathbf{s}}

\newtheorem{theorem}{Theorem}[section]
\newtheorem{lemma}{Lemma}[section]
\newtheorem{proposition}{Proposition}[section]

\begin{document}

\def\spacingset#1{\renewcommand{\baselinestretch}%
{#1}\small\normalsize} \spacingset{1}


\if1\blind
{
  \title{\bf Objective Bayesian analysis for spatial Student-t regression models}
  \author{Jose A. Ordoñez\\
    Departamento de Estatística, Universidade Estadual de Campinas\\
    Marcos O. Prates\thanks{Correspondece to: Marcos O. Prates, Statistics Department, Universidade Federal de Minas Gerais, Av. Ant\^onio Carlos,  6627 - Pr\'edio do ICEx - sala 4077, Belo Horizonte, Minas Gerais, Brasil. Email: marcosop@est.ufmg.br} \\
    Departamento de Estatística, Universidade Federal de Minas Gerais\\
    Larissa A. Matos\\
    Departamento de Estatística, Universidade Estadual de Campinas\\
    Victor H. Lachos \\
    Department of Statistics, University of Connecticut\\
    }
    
    \maketitle
} \fi

\if0\blind
{
  \bigskip
  \bigskip
  \bigskip
  \begin{center}
    {\LARGE\bf Objective  Bayesian  prior  for  geostatistical Student-t models}
\end{center}
  \medskip
} \fi
\bigskip
\begin{abstract}
The choice of the prior distribution is a key aspect of Bayesian analysis. For the spatial regression setting a subjective prior choice for the parameters may not be trivial, from this perspective, using the objective Bayesian analysis framework a reference is introduced for the spatial Student-t regression model with unknown degrees of freedom. The spatial Student-t regression model poses two main challenges when eliciting priors: one for the spatial dependence parameter and the other one for the degrees of freedom. It is well-known that the propriety of the posterior distribution over objective priors is not always guaranteed, whereas the use of proper prior distributions may dominate and bias the posterior analysis. In this paper, we show the conditions under which our proposed reference prior yield to a proper posterior distribution. Simulation studies are used in order to evaluate the performance of the reference prior to a commonly used vague proper prior.
\end{abstract}

\noindent%
\vfill

\newpage
\spacingset{1.45} 

\section{Introduction}

Geostatistical data modeling \citep{cressie:1993} has now virtually permeated all areas of epidemiology, hydrology, agriculture, environmental science, demographic studies, just to name a few. Here, the prime objective is to account for the spatial correlation among observations collected at various locations, and also to predict the values of interest for non-sampled sites. In this paper we will focus in a fully Bayesian approach to analyze spatial data, whose main advantage is that parameter uncertainly is fully accounted for when performing prediction and inference, even in small samples \citep{berger}. However, elicitation of priors for correlation parameter in a Gaussian processes is a non trivial task \citep{kennedy2001bayesian}.

The problem of inference and prediction for spatial data with Gaussian processes using objective priors has received attention in the recent literature. It started with \citet{berger} that develop an exact non-informative prior for unknown parameters of Gaussian random fields by using exact marginalization in the reference prior algorithm \citep[See,][]{bernardo1992}. Further, \citet{paulo2005default} and  \citet{ren2013objective} generalized the previous results to an arbitrary number of parameters in the correlation parameter. After the precursor proposal of \citet{de2007objective} that allow the inclusion of measurement error for reference prior elicitation, other extension from this perspective were proposed in the literature, see for instance, \cite{ren2012objective} and \cite{kazianka2012objective}.

In the context of the Student-t distribution, \citet{zellner1976} was the first to present a Bayesian and non-Bayesian analysis of a linear multiple regression model with Student-t errors assuming a scalar dispersion matrix and known degrees of freedom. An interesting result of this paper is that 
inferences about the scale parameter of the multivariate-t distribution can be made using an F-distribution rather than the usual $\chi^2$(or inverted $\chi^2$) distribution. Later, \citet{fonseca2008} developed an objective Bayesian analyses based on the Jeffreys-rule prior and on the independence Jeffreys prior for linear regression models with independent Student-t errors and unknown degrees of freedom. This procedure allowed a non-subjective statistical analysis with adaptive robustness to outliers and with full account of the uncertainty.  \citet{branco2013objective} introduced an objective prior for the shape parameter using the skew-t distribution proposed by \cite{azzalini2003distributions}. More recently, \cite{villa2014} constructed an objective prior for the degrees of freedom of the univariate Student-t distribution when this parameter is taken to be discrete. 

Even though some solutions have been proposed in the literature to deal with the problem of
objective prior under the Student-t distribution, to the best of our knowledge there are not studies conducting objective Bayesian analyzes under the Student-t spatial regression model. Following \citet{berger}, we introduce a reference prior based on exact marginalization and we derive the conditions that it yields a valid posterior distribution. Moreover, the independence Jeffrey and the Jeffrey-rule priors are derived and analyzed. As in \citet{berger}, we show that the Jeffreys priors suffers many drawbacks while the proposed reference prior produces more accurate estimates with good frequentist properties.  

The paper is organized as follows. In Section~\ref{s:pre}, we describe the Student-t spatial regression model as well as the family of covariance functions that will be considered. In Section~\ref{s:ref}, a general form of improper priors are presented and the reference prior is provided with the conditions of its validity.  In Section~\ref{s:ms}, model selection criteria are presented in order to evaluate the competing Bayesian models. In Section~\ref{s:sim}, a simulation study is performed to assess the frequentist properties of the Bayesian estimates under different priors. Finally, a brief discussion is presented in Section~\ref{s:conc}.

\section{The Student-t Spatial Regression Model}
\label{s:pre}
Let $Y(\sa)$ denote the response over location $\sa \in D_{\sa}$, where $D_{\sa}$ is a continuous spatial domain in ${\rm I\!R}^2$. We assume that, the observed data $\mathbf{y}(\sa) = (y(\sa_1), \ldots, y(\sa_n))^\top$ is a single realization of a  Student-t stochastic process, $Y(\sa) \equiv {y(\sa): \sa \in D_{\sa}}$ \citep{palacios2006non,bevilacquanon}. Thus, if $\mathbf{Y}$ follow a multivariate Student-t distribution with location vector $\bmu$, scale matrix $\bSigma$ and $\nu$ degrees of freedom, $\Y \sim t_n(\X\bbeta,\bSigma,\nu)$, the Student-t spatial regression (T-SR) model can be represented as
\begin{eqnarray*}
\Y=\X\bbeta+\bepsilon, 
\end{eqnarray*}
where 
$\X$ is a $n \times p$ non-stochastic matrix of full rank with the $i$th row $\xp_i^{\top}=(x_{i1},\ldots,x_{ip})^{\top}=(x_1(\sa_i),\ldots,x_p(\sa_i))^{\top}$, $\bbeta=(\beta_1,\ldots,\beta_p)^{\top}$, and $\bepsilon=(\epsilon(\sa_1),\ldots,\epsilon(\sa_n))^{\top}$ with  $E[\bepsilon]=\mathbf{0}$.

Equivalently, the model can be written as 
\begin{eqnarray}\label{model}
Y(\sa_i)= \bmu(\sa_i)+\epsilon(\sa_i), \quad i= 1,\ldots, n,
\end{eqnarray}
where  $\bmu(\sa_i)=\sum_{j}x_j(\sa_i)\beta_j$ is the mean of the stochastic process for $j=1,\ldots,p$.
Therefore, a realization of a Student-t process can be represented by setting $\bepsilon \sim t_n(\mathbf{0},\bSigma,\nu)$ with a valid covariance function for the scale matrix $\bSigma$. We concentrate on a particular parametric class of covariance functions such that the scale matrix is given by $$\bSigma=[C(\sa_i,\sa_j)]=\sigma^2\R(\phi) + \tau \mathbf{I},$$ 
where in standard geostatistical terms: $\sigma^2$ is the sill; $\tau$ is the nugget effect; $\phi$ determines the range of the spatial process; $\R(\phi)$ is an $n\times n$ correlation matrix; and $\mathbf{I}$ is the identity matrix. We assume that $\R(\phi)$ is an isotropic correlation matrix and depends only the Euclidean distance $d_{ij} = ||\sa_i - \sa_j||$ between the points $\sa_i$ and $\sa_j$. Thus, the likelihood function of the model parameters $(\bbeta,\sigma^2,\bphi,\tau,\nu)$, based on the observed data $\yp$, is given by

\begin{equation}\label{like}
L(\beta,\sigma^2,\phi,\tau, \nu|\yp)= \frac{\bGamma\left(\frac{\nu+n}{2}\right)\nu^{\frac{\nu}{2}}}{\bGamma\left(\frac{\nu}{2}\right)\pi^{\frac{n}{2}}|\bSigma|^{\frac{1}{2}}}\left \lbrace \nu+ (\yp-\X\bbeta)^{\top}\bSigma^{-1}(\yp-\X\bbeta)\right \rbrace ^{-\frac{\nu+n}{2}},
\end{equation}
where $|\mathbf{A}|$ denotes the determinant of the matrix $\mathbf{A}$. In this work, we consider four general families of isotropic correlation functions in ${\rm I\!R}^2$, says, spherical, Cauchy, power exponential and Mat\'ern correlations functions. The spherical family have one difference which distinguish it from all the other families, this family present a finite range i.e., $\R(\phi)= 0$ for $d_{ij} > \phi$. The power exponential and Mat\'ern have the traditional exponential and Gaussian covariance functions as special cases. For further details, see \citet{banerjee2014hierarchical}.

\section{The Reference Prior}
\label{s:ref}

\citet{palacios2006non} discuss that the derivation of a reference prior for non-Gaussian processes is not trivial. In this section, we introduce a reference prior for the T-SR model defined in ~\eqref{model} without a nugget effect, i.e., $\bSigma=\sigma^2\R(\phi)$. 
We obtain $\pi(\phi,\nu)$ through the marginal model defined via integrated likelihood.

\subsection{Prior density for $(\bbeta,\sigma^2,\phi,\nu)$}

For $(\bbeta,\sigma^2,\phi,\nu) \in \Omega= {\rm I\!R}^p \times (0,\infty)\times (0,\infty)\times (0,\infty)$, consider the family of improper priors of the form:
\begin{equation}\label{priorform}
\pi(\bbeta,\sigma^2,\phi,\nu) \propto \frac{\pi(\phi,\nu)}{(\sigma^2)^{a}},
\end{equation}
for different choices of $\pi(\phi,\nu)$ and $a$. Selection of the prior distribution for $\phi$ and $\nu$ is not straightforward. Assuming an independence structure $\pi(\phi,\nu)= \pi_1(\phi) \times \pi_2(\nu)$, one alternative could be to select improper priors for these two parameters, nevertheless, it is necessary to be careful, since it is obligatory to show that such selection produces a proper posterior distribution. For $\phi$, the use of truncation over the parameter space or vague proper priors are alternatives to overcome the improper posterior distribution problem, however, in both cases inferences are often highly dependent on the bounds used or on the hypeparamters selected for the vague distribution \citep{berger}. And, for $\nu$ even when the parameter space is restricted, the maximum likelihood estimator may not exist with positive probability \citep{fonseca2008}. Other choices of priors can be found in \citet{Geweke1993,de2007objective,ren2012objective,kazianka2012objective,ren2013objective,fonseca2008,branco2013objective}, but none of this authors consider a Student-t spatial regression framework.

The expression for $\pi(\phi, \nu| \yp)$ based on an arbitrary prior $\pi(\phi,\nu)$ is presented in the following proposition.
\begin{proposition} \label{prop2}
 For $a< \nicefrac{\nu}{2}+1$ and different choices for $\pi(\phi,\nu)$, the posterior density $\pi(\phi, \nu| \yp)$ can be written as
\begin{equation}\label{posteriorphi}
\pi(\phi,\nu| \yp) \propto A(\nu)\left|\R\right|^{-\frac{1}{2}}\left|V_{\hat{\bbeta}}\right|^{\frac{1}{2}}\left \lbrace (n-p) S^2\right  \rbrace^{-(\frac{n-p}{2}+a-1)}\pi(\phi,\nu),
\end{equation}
where $A(\nu)=\nu^{-(1-a)}{\Gamma(\frac{\nu}{2}-a+1)}/{\Gamma(\frac{\nu}{2})}$, $S^2=\displaystyle\nicefrac{\zp^{\top}\R^{-1}\zp}{(n-p)}$, $\zp=(\yp-\X\hat{\bbeta})$ and $V_{\hat{\bbeta}}=(X^{\top}\R^{-1}X)^{-1}$. Let $\hat{\bbeta}=(\X^{\top}\mathbf{R}^{-1}X)^{-1}\X^{\top}\mathbf{R}^{-1}\y$ be the generalized least square estimator of $\bbeta$. So, given $\pi(\phi,\nu)$, to guarantee propriety of the posterior density $\pi(\bbeta,\sigma^2,\phi,\nu|\yp)$ and the existence of the first two moments of the Student-t distribution, we have to ensure that $\forall \varepsilon > 0$,
$$\int_{0}^{+\infty} \int_{4+\varepsilon}^{+\infty}  A(\nu)\left|\R\right|^{-\frac{1}{2}}\left|V_{\hat{\bbeta}}\right|^{\frac{1}{2}}\left \lbrace(n-p)S^2\right \rbrace^{-(\frac{n-p}{2}+a-1)}\pi(\phi,\nu)d \nu d \phi < +\infty.$$
\end{proposition}

\subsection{Proposal of $\pi(\phi,\nu)$}\label{sub32}

Let $\theta=(\bbeta,\theta^*)$, with  $\theta^{*}=(\sigma^2,\phi,\nu)$ and $L(\theta|\y)$ the sampling distribution defined in \eqref{like}. For the reference prior, $\theta^*$ is the parameter of interest and we assume that $\bbeta$ is a nuisance parameter. Now, factorizing the prior distribution $\pi(\bbeta,\theta^{*})=\pi(\bbeta| \theta^{*})\pi(\theta^{*})$ and choosing $\pi(\bbeta|\theta^{*})=1$ as this is the reference prior in Equation~\eqref{like}, we have that 
\begin{equation}
\label{eqL1}
L_1(\theta^{*})=\int_{{\rm I\!R}^{p}}L(\theta|\yp)\pi(\bbeta|\theta^{*})d\bbeta \propto  \frac{\Gamma(\frac{\nu_1}{2})}{\Gamma(\frac{\nu}{2})}\frac{(\sigma^2)^{-\frac{n-p}{2}}}{\left\lbrace\nu + (n-p) S^{*}\right \rbrace^{\frac{\nu_1}{2}}}|V_{\hat{\bbeta}}|^{\frac{1}{2}}|\R|^{-\frac{1}{2}}\nu^{\frac{\nu}{2}},
\end{equation}
where $S^{*}=\nicefrac{S^2}{\sigma^2}$ and $\nu_1= (n-p)+\nu$. It is possible to show that this expression converge to the normal case when $\nu \to +\infty$. Using the prior reference method \citep{bernardo1992}, it is necessary to calculate $E\left \lbrace \Z^{\top}\bSigma^{-1}\frac{\partial \bSigma}{\partial \phi}\bSigma^{-1}\Z|S^{*}\right \rbrace$ and $E\left \lbrace \left(\Z^{\top}\bSigma^{-1}\frac{\partial \bSigma}{\partial \phi}\bSigma^{-1}\Z\right)^2|S^{*}\right \rbrace$. Unfortunately, these conditional expectations have no analytical form for the Student-t case. One possible solution is to numerically compute these expressions by using Monte Carlo approximation which will demand a high computational cost making inference infeasible. 
For this reason, we suggest the use of the marginal expectations $E\left \lbrace \Z^{\top}\bSigma^{-1}\frac{\partial \bSigma}{\partial \phi}\bSigma^{-1}\Z\right\rbrace $ and $E\left \lbrace \left(\Z^{\top}\bSigma^{-1}\frac{\partial \bSigma}{\partial \phi}\bSigma^{-1}\Z\right)^2 \right \rbrace$ in our prior proposal. This suggestion may result in a improper prior (Theorem \ref{prop3}), but lead to a proper posterior distribution (Theorem \ref{prop5}).

\begin{theorem}\label{prop3} Under the T-SR model defined in (\ref{model}), for $\phi>0$ and $\nu> 4 + \varepsilon$, for any $\varepsilon>0$, the prior distribution obtained through the reference prior method is of the form (\ref{priorform}), with $a=1$ and 
\begin{equation}\label{proposalprior}
\pi(\phi,\nu)\propto \left(BCD+16\left(B_{11}C_{11}B_{12}-BC_{11}^2\right)- 8B_{12}^2C- \frac{1}{2}B_{11}^2D\right)^{\frac{1}{2}},
\end{equation}
where
\begin{eqnarray*}
B&=&\frac{\nu (n-p)}{\tau(\nu)},\,\, C=\left(\frac{2 (n-p)}{\tau(\nu)\nu}+1\right)\mathbf{A}-\frac{\nu+2}{\nu-2}tr^2\left[\Phi \right],\\
D&=& -\left(\frac{2(n-p)}{\nu} \frac{\tau(\nu)+2}{\tau(\nu)(\tau(\nu)-2)}  +\delta_1(\nu)\right),\,\,B_{11}=-\frac{2\nu (n-p)}{(\nu-2)\tau(\nu)}tr\left[\Phi\right],\\
B_{12}&=&-\frac{n-p}{(\tau(\nu)-2)\tau(\nu)}\,\, \mbox{\it and}\,\,C_{11}=\frac{n-p}{(\nu-2)(\tau(\nu)-2)\tau(\nu)}tr \left[\Phi\right],
\end{eqnarray*}
with $\tau(\nu)= n-p+\nu+2$,\,$\Phi=\frac{\partial \R}{\partial \phi}\R^{-1}\mathbf{P},$\,\,
$\mathbf{A}=\displaystyle\frac{\nu^2}{(\nu-2)(\nu-4)}\left(2tr\left[\Phi^2\right] + tr^2\left[\Phi \right]\right)$, $\mathbf{P}=I-\X(\X^{\top}\R^{-1}\X)^{-1}\X^{\top}\R^{-1}$\,\,
and
$\delta_1(\nu)=\Psi_1(\displaystyle\nicefrac{\nu+ n-p}{2})-\Psi_1(\nicefrac{\nu}{2})$, where $\Psi_1(.)$ denotes the trigamma function.
\end{theorem} 
 
A short proof of this theorem can be found in Appendix~\ref{ap1}. The following Lemma provides conditions to show the results of Theorem~\ref{prop5}.
\begin{lemma} \label{prop4}
For any $\varepsilon > 0$, we have that:
\begin{itemize}
\item If $\mathbf{1}$ is not a column of $\X$. Then,
\begin{enumerate}[(a)]

\item as $\nu \to +\infty$ and $\phi \to +\infty$,
$$\pi(\phi,\nu) = O\left(\nu^{-\frac{3}{2}}\frac{d}{d \phi}\log\psi(\phi)\right);$$

\item as $\nu \to 4+\varepsilon$ and $\phi \to +\infty$,
$$\pi(\phi,\nu) = O\left(\frac{d }{d \phi}\log\psi(\phi)\right).$$
\end{enumerate}

\item If $\mathbf{1}$ is a column of $\X$. Then,

\begin{enumerate}[(a)]
\item as $\nu \to +\infty$ and $\phi \to +\infty$,
$$\pi(\phi,\nu) = O\left(\nu^{-\frac{3}{2}}\frac{\omega(\phi)}{\psi(\phi)}\frac{d }{d \phi}\log \left[\frac{\omega(\phi)}{\psi(\phi)}\right]\right);$$

\item as $\nu \to 4+\varepsilon$ and $\phi \to +\infty$,

$$\pi(\phi,\nu) = O\left(\frac{\omega(\phi)}{\psi(\phi)}\frac{d }{d \phi}\log \left[\frac{\omega(\phi)}{\psi(\phi)}\right]\right),$$
\end{enumerate}
\end{itemize}
where $\psi(\phi)=\nu(\phi)$ and $\omega(\phi)$ are as defined in \citet{berger} for each correlation matrix considered in Section~\ref{s:pre}. 
\end{lemma}

See Appendix~\ref{ap2} for the proof.  The next Theorem~\ref{prop5} shows the conditions under which the reference prior introduced in Theorem~\ref{prop3} is proper.

\begin{theorem}\label{prop5}
For any of the families of correlation function considered in Section~\ref{s:pre} and under the T-SR model~\eqref{model}, the posterior distribution of $(\bbeta, \sigma, \phi, \nu)$ is  proper if the conditions of in Table~\ref{tabla1} are satisfied for the hyperparameter $a$.

\begin{table}[htbp]
  \centering
  \caption{Conditions to guarantee the propriety of the posterior distribution using the proposal prior.}
    \begin{tabular}{lll}
   \hline
    Correlation Family & 1 is a column of X & 1 is not a column of X \\
   \hline
    Spherical  & $-1<a<\frac{\nu}{2}+1$  & $\frac{1}{2}<a<\frac{\nu}{2}+1$ \\
    Power exponential & $0 < a< \frac{\nu}{2}+1$  & $\frac{1}{2}<a < \frac{\nu}{2}+1$  \\
    Cauchy & $0<a < \frac{\nu}{2}+1$ & $\frac{1}{2}<a  < \frac{\nu}{2}+1$ \\
    Matern &       &  \\
    $\kappa <1$ & $2- \kappa^{-1} <a < \frac{\nu}{2}+1 $ & $\frac{1}{2}< a<\frac{\nu}{2}+1$ \\
    $\kappa \geq 1$ & $\kappa^{-1}<a<\frac{\nu}{2}+1 $ & $\frac{1}{2}< a <\frac{\nu}{2}+1$ \\
    \hline
    \end{tabular} \label{tabla1}  
    
\end{table}%

\end{theorem}

\textit{Proof:} To guarantee propriety we have to ensure that Propostion~\ref{prop2} is satisfied. Let $\varepsilon > 0$ and
\begin{equation}\label{l1}
    L^{*}(\phi, \yp)= \left|\R\right|^{-\frac{1}{2}}\left|V_{\hat{\bbeta}}\right|^{\frac{1}{2}}(\nu S^2)^{-(\frac{\nu}{2}+a-1)}.
\end{equation}
One can show that $A(\nu)$ has a constant behavior when $\nu \to +\infty$. Also, we have that $\nu^{-\nicefrac{3}{2}}$ is integrable at $(4 + \varepsilon,+\infty)$. Using both facts and combining with Lemma~\ref{prop4}, we have that when $\mathbf{1}$ is not a column of $\X$, then
\begin{eqnarray} \nonumber
\int_{0}^{+\infty} \int_{4+\varepsilon}^{+\infty} && A(\nu)L^{*}(\phi,\yp)\pi(\nu,\phi)d \nu d\phi \\ \nonumber
&\leq& C_0\int_{0}^{+\infty}L^{*}(\phi,\yp)\frac{d}{d \phi} \log\psi(\phi)d\phi \int_{4+\varepsilon}^{+\infty}\nu^{-\frac{3}{2}}d \nu \\ \nonumber
&\leq& C_1 \int_{0}^{+\infty}L^{*}(\phi,\yp)\frac{d }{d \phi}\log\psi(\phi)d\phi < +\infty;
\end{eqnarray}
and, when $\mathbf{1}$ is a column of $\X$, then
\begin{eqnarray} \nonumber
\int_{0}^{+\infty} \int_{4+\varepsilon}^{+\infty} && A(\nu)L^{*}(\phi,\yp)\pi(\nu,\phi)d \nu d\phi \\ \nonumber
&\leq& C_2\int_{0}^{+\infty}L^{*}(\phi,\yp)\frac{\omega(\phi)}{\psi(\phi)}\frac{d }{d \phi}\log \left[\frac{\omega(\phi)}{\psi(\phi)}\right]d\phi  \int_{4+\varepsilon}^{+\infty}A(\nu)\nu^{-\frac{3}{2}}d \nu  \\ \nonumber
& \leq & C_3 \int_{0}^{+\infty}L^{*}(\phi,\yp)\frac{\omega(\phi)}{\psi(\phi)}\frac{d }{d \phi}\log\left[\frac{\omega(\phi)}{\psi(\phi)}\right]d\phi < +\infty.
\end{eqnarray}
with $C_k, \,\,k=0,1,2,3$ being constants, both integrals are finite as a consequence of Lemma 1 and Theorem 4 given in \citet[][pg. 1366]{berger}. 

\section{Model selection}
\label{s:ms}

Let us start by setting up the model selection as an hypothesis testing problem \citep{banerjee2014hierarchical,berger}. Thus, replace the usual hypotheses by a candidate parametric model, say $m_k$, having respective parameter vectors $\theta_{m_k}$.
Under the prior density proposal in Equation~\eqref{priorform} we compute the marginal density for a model $m_k$ as
\begin{eqnarray} \nonumber
m_k(\y) &=& \int L(\theta_{m_k})\pi(\theta_{m_k})d \theta_{m_k}\\ \nonumber &=&\int_{0}^{+\infty}\int_{4+\varepsilon}^{+\infty} A(\nu)\left|\R\right|^{-\frac{1}{2}}\left|V_{\hat{\bbeta}}\right|^{\frac{1}{2}}\left \lbrace(n-p) S^2\right \rbrace^{-\left(\frac{n-p}{2}+a-1\right)}\pi(\phi,\nu)d \nu d\phi.
\end{eqnarray}

To compare $q$ different models $m_{k}$, with $k=1,\ldots,q$, we assign equal prior probabilities to the models. Therefore, the resulting posterior probability for the $k$-th model is defined by
\begin{equation*}
p(m_k \mid \y) = \frac{m_{k}(\y)}{\sum_{j=1}^{q} m_j(\y)}.
\end{equation*}
Under this criterion a model with the largest posterior probability is preferable. Another possibility to perform model selection is to choose the model with best prediction power. Suppose that $n_0$ locations are separated as a validation set. The mean square prediction error (MSPE) of the $k$-th model is defined by
$$
\text{MSPE}_k = \frac{\sum_{i=1}^{n_0} \left \lbrace {Y(s_i) - \hat Y^{(k)}(s_i)}\right \rbrace^2}{n_0},
$$
where $\hat Y^{(k)}(s_i)$ is the predicted response for observation in location $s_i$ under the $k$-th model. Thus, the model that minimizes the MSPE is the most suitable under this criterion.

\section{Simulation Study}
\label{s:sim}

The study of the frequentist properties of Bayesian inference is of interest to assess and understand the properties of non-informative or default priors \citep[e.g.,][]{stein1985,berger2006case,de2007objective,kazianka2012objective,branco2013objective,he2020objective}. Therefore, a simulation study is performed to assess the performance of the proposal method and compare it with a vague prior ({\it vague}) of the form \eqref{priorform}.

One of the most used and flexible isotropic correlation function is the Mat\'{e}rn. The Mat\'{e}rn family is defined as
\begin{equation}\label{Rphi}
\R(\phi)=\left\{\begin{array}{ll}
\displaystyle \frac{1}{2^{\kappa-1}\Gamma{(\kappa)}}\left(\frac{d_{ij}}{\phi}\right)^{\kappa}K_{\kappa}(d_{ij}/\phi),& d_{ij}>0, \\
1, & d_{ij}=0,
\end{array}\right.
\end{equation}
where $\phi>0$; $K_{\kappa}(u)=\frac{1}{2}\int_{0}^{\infty}x^{\kappa-1}e^{-\nicefrac{1}{2}u(x + x^{-1})}$ is the modified Bessel function of the third kind of order $\kappa$ \citep[see][]{gradshtejn1965table}, with $\kappa > 0$ fixed. When $\kappa\rightarrow\infty$  and $\kappa=0.5$, the Gaussian and exponential correlations respectively can be obtained from (\ref{Rphi}) \citep[see,][]{diggle2007springer}.

We propose two T-SR models with coordinates $\sa=(\mathbf{x}_1, \mathbf{x}_2)$ and $\R$ belonging to the Mat\'{e}rn family with $\kappa =0.5$ (exponential correlation structure), $\sigma^2=0.8$ and $\phi=2$ to study the proposed  priors. The first one ({\bf Scenario 1}) is given by, 
\begin{equation}
\label{modela}
y(s_i) = 10 + \varepsilon_{s_i}, \quad i, \ldots,n, 
\end{equation}
with $\mathbf{\varepsilon} \sim t_n(\mathbf{0}, \sigma^2\R,\nu=5)$. And, to illustrate the beyond a trivial intercept model, the second T-SR model ({\bf Scenario 2}) is given by
\begin{equation}
y(s_i) = 0 -2.2x_{1i} + 0.5x_{2i}  + 1.7x_{1i}^2 + 2.4x_{2i}^2 + 3.5x_{1i}x_{2i} + \varepsilon_{s_i}.
\label{modelb}
\end{equation}
A total of $K=500$ Monte Carlo simulations were generated for each scenario, the coordinates $\sa$ were sampled at $n=100$ locations of a regular lattice in $ D_{\sa}=[0,10] \times [0,10]$. 

For the vague proper prior, we consider $a=2.1$ and $\pi(\phi,\nu)=\pi(\phi)\times \pi(\nu)$ with $\pi(\phi)= U(0.1,4.72)$ and $\pi(\nu)= Texp(\lambda; \nu \in {\rm A})$, where $\pi(\lambda) = U(0.01,0.25)$, ${\rm A}=[4.1,+\infty)$ and $Texp(\lambda, \nu \in {\rm A})$ denote the truncated exponential distribution. The distribution of $\lambda$ is such that it allows the mean of the prior of $\nu$ to vary from $4$ to $100$. The prior exponential prior of $\nu$ is truncated above $4.1$ to guarantee the existence of the Student-t process and the prior of $\phi$ allow that the distance such that the empirical range, corr$(s_i, s_j) < 0.05$, varies from $0.30$ to $14$ (which is the minimum and maximum distance between the locations, respectively). 

For the two scenarios, we compute the empirical equal-tailed 95$\%$ credible interval for all parameters, based on the two priors. We also compute the coverage probability for each parameter as the number of simulations in such the parameter is inside the credible limits, and the expected log length of each credible interval as the mean of the logarithm of the difference (Log-length) between the upper and lower credible limits for each simulation. The bias of each parameter was estimated as $\text{Bias}_j = \sum_{k=1}^K (\hat \theta^{k}_j - \theta_j)/K$, where $\theta_j$ is the true parameter value and $\hat \theta^{k}_j$ is the median posterior estimate for the $j$-th parameter in the $k$-th Monte Carlo simulation.
\begin{table}[htbp]
  \centering
  \caption{Simulation Results under {\bf Scenario 1}.}
  \label{t:sim1}
  \begin{scriptsize}
    \begin{tabular}{llcccc}
    \hline
     \multicolumn{1}{l}{} & \multicolumn{1}{c}{Prior} & \multicolumn{1}{c}{$\beta_0$} & \multicolumn{1}{c}{$\sigma^2$} & \multicolumn{1}{c}{$\phi$} & \multicolumn{1}{c}{$\nu$} \\
  \hline
   \multirow{2}[2]{*}{Bias} & {\it reference} & -0.037 (0.41) & 0.465 (1.31) & 0.140 (0.58) & 0.795 (0.29) \\
      & {\it vague}& -0.038 (0.42) & 4.026 (9.62) & 0.326 (0.89) & 3.662 (1.69) \\ 
 \hline
\multirow{2}[2]{*}{Log length} & {\it reference}  & 0.08 & 1.034 & 0.974 & 2.442 \\ 
      & {\it vague}& 0.062 & 0.733 & 0.8161 &  2.360 \\ 
\hline
\multirow{2}[2]{*}{C.P} & {\it reference} & 0.945 & 0.948 & 0.988 & 1.000 \\
        & {\it vague} & 0.970 & 0.5831 & 0.878 & 0.907 \\
      \hline
    \end{tabular}%
  \end{scriptsize}
 \label{taba}%
\end{table}%

Table~\ref{t:sim1} shows the results under {\bf Scenario 1}. As we can see, there is almost no Bias for $\beta_0$ with compatible standard deviations for all priors. However, when looking over the hyperparamteres the {\it vague} prior have larger bias for $\sigma^2$ and $\nu$. The {\it reference} prior is the one that better achieve the nominal coverage of 95\% for all parameters.
\begin{table}[htbp]
  \centering
  \caption{Simulation results under {\bf Scenario 2}.}
  \label{t:sim2}
\resizebox{.98\textwidth}{!}{%
    \begin{tabular}{llccccccccc}
    \hline
      \multicolumn{1}{l}{} & \multicolumn{1}{c}{Prior} & $\beta_0$ & $\beta_1$ & $\beta_2$ & $\beta_3$ & $\beta_4$ & $\beta_5$ &$\sigma^2$ & $\phi$   & $\nu$  \\
      \hline
\multirow{2}[0]{*}{Bias} & {\it reference} & 0.016 (1.09) & -0.04 (0.33) & -0.04 (0.33) & 0.002 (0.02) & 0.002 (0.02) & 0.002 (0.02) & 0.12 (0.87) & -0.58 (0.38) & 0.7393 (0.52) \\ 

      & {\it vague} & 0.04 (1.11) & -0.03 (0.35) & -0.04 (0.34) & 0.003 (0.03) & 0.002 (0.03) & 0.002 (0.03) & 4.37 (8.3) & 0.56 (0.78) & 3.8013 (1.91)\\ 
\hline
\multirow{2}[0]{*}{Log length} & {\it reference} & 0.4526 & 0.092 & 0.111 & 0.9946 & 1.186 & 0.3552 & 0.3270 & 0.3472 & 0.930 \\
      & {\it vague} & 0.623 & 0.08 & 0.0582 & 0.8404 & 0.9681 & 0.2992 & 0.543 & 0.2965 & 0.7824 \\
\hline
\multirow{2}[0]{*}{C.P} & {\it reference} & 0.958 & 0.916 & 0.904 & 0.936 & 0.894 & 0.966 & 0.924 & 0.952 & 1 \\
      & {\it vague} & 0.942 & 0.936  & 0.94 & 0.945 & 0.921 & 0.873 & 0.623     & 0.615  & 0.99 \\
\hline
\end{tabular}}%
 \label{tabb}%
\end{table}%

Table~\ref{t:sim2} show the results under {\bf Scenario 2}. For all priors the estimates of $\beta$'s seems to perform in a similar manner. The {\it reference} and {\it vague} priors provide closer coverage probability to the nominal value. Focusing in the hyperparameters $\theta^*$, 
it is clear again that the proposed reference prior is the one that provides the best results. Overall, it has smaller bias than its competitor with and adequate credible interval length, since it is the only prior that provides coverages close to the nominal value for all parameters. 

\section{Discussion}
\label{s:conc}

In this paper we propose and recommend a reference prior for the spatial Student-t regression. For the proposed prior, the conditions under which it yields a proper posterior distribution were presented and discussed.

We show through simulations that the {\it reference} prior presents better performance than the {\it vague} prior . It shows small estimation bias and adequate frequentist coverage for all parameters. The \texttt{OBASpatial} \texttt{R} package is available at \texttt{CRAN} for download and allow practitioners to fit the proposed model with the different priors introduced in the manuscript. 

As further studies, the inclusion calculation of Jeffrey's priors for the Student-t spatial regression is of interest Different generalizations of the family of correlation functions are also of interes.

\section*{Funding}\label{ack}
The research of Jose A. Ordoñez is funded by CAPES. Marcos O. Prates would like to acknowledge CNPq grants 436948/2018-4 and 307547/2018-4 and FAPEMIG grant PPM-00532-16 for partial financial support.

\section*{Appendix}\label{app}

\subsection{Proof of Theorem~\ref{prop3}}
\label{ap1}

Let $\ell_1(\theta^{*})=\log(L_1(\theta^{*}))$, with $L_1(\theta^{*})$ as defined in (\ref{eqL1}). By using  the marginal expectations $E\left \lbrace \Z^{\top}\bSigma^{-1}\frac{\partial \bSigma}{\partial \phi}\bSigma^{-1}\Z\right\rbrace $ and $E\left \lbrace \left(\Z^{\top}\bSigma^{-1}\frac{\partial \bSigma}{\partial \phi}\bSigma^{-1}\Z\right)^2 \right \rbrace$  as suggested in subsection~\ref{sub32}, we have that 
\begin{equation*}
I_{1}(\theta^{*})=
\begin{bmatrix}
 I_{1(\sigma^2)^2} & I_{1\sigma^2 \phi} & I_{1\sigma^2 \nu}\\
I_{1\sigma^2 \phi} & I_{1(\phi)^2} & I_{1\phi \nu}\\
I_{1\sigma^2 \nu} & I_{1\phi \nu} & I_{1\nu^2}
\end{bmatrix}.
\end{equation*}

We also have 
\begingroup
\allowdisplaybreaks
\begin{eqnarray*}
I_{1(\sigma^2)^2} &  =&\frac{1}{2\sigma^4}\frac{\nu n_p}{\tau(\nu)}=\frac{1}{2\sigma^4}B, \\
I_{1\phi^2} & = &\frac{1}{4}\left \lbrace \left(\frac{2 n_p}{\tau(\nu)\nu}+1\right)\mathbf{A}-\left(\frac{\nu+2}{\nu-2}\right)tr^2\left[\Phi \right] \right \rbrace=\frac{1}{4}C,\\
I_{1\nu^{2}} &  =&\frac{1}{4}\left\lbrace -\frac{2n_p}{\nu} \left(\frac{\tau(\nu)+2}{\tau(\nu)(\tau(\nu)-2)}  \right)-\delta_1(\nu)\right\rbrace= \frac{1}{4}D,\\
I_{ 1\sigma^2 \phi} & =& -\frac{1}{4\sigma^2}\left ( \frac{2\nu n_p}{(\nu-2)\tau(\nu)}tr\left[\Phi\right]\right)= \frac{1}{4\sigma^2}B_{11} ,\\
I_{ 1\sigma^2 \nu} & =&-\frac{n-p}{\sigma^2(\tau(\nu)-2)\tau(\nu)}=\frac{1}{\sigma^2}B_{12},\\
I_{ 1\phi \nu} & =& \frac{n-p}{(\nu-2)(\tau(\nu)-2)\tau(\nu)}tr \left[\Phi\right]= C_{11},
\end{eqnarray*}
\endgroup
where   $\mathbf{A}=\frac{\nu^2}{(\nu-2)(\nu-4)}\left\lbrace 2tr\left[\Phi^2\right] + tr^2\left[\Phi \right]\right\rbrace$, \,$\tau(\nu)= n-p+\nu+2$,\,$\Phi=\frac{\partial \R}{\partial \phi}\R^{-1}\mathbf{P}$, 
and $\delta_1(\nu)=\Psi_1(\nicefrac{\nu+ n-p}{2})-\Psi_1(\nicefrac{\nu}{2})$, whit $\Psi_1(.)$ denoting the trigamma function. Therefore,

\begin{equation*}
 I_1(\theta^{*})=\begin{bmatrix}
 \frac{1}{2\sigma^4}B & \frac{1}{4\sigma^2}B_{11} & \frac{1}{\sigma^2}B_{12},\\
 \frac{1}{4\sigma^2}B_{11} & \frac{1}{4}C & C_{11},\\
 \frac{1}{\sigma^2}B_{12} & C_{11} & \frac{1}{4}D
 \end{bmatrix}
 \end{equation*}
and 
 \begin{equation*}
\pi(\theta)= |I_1(\theta^{*})|^{\frac{1}{2}}  \propto \frac{1}{\sigma^2}\left(BCD+16B_{11}C_{11}B_{12}- 8B_{12}^2C- 16BC_{11}^2- \frac{1}{2}B_{11}^2D\right)^{\frac{1}{2}}.
 \end{equation*}
 
\subsection{Proof of Lemma~\ref{prop4}}\label{ap2}

Previously, in Proposition \ref{prop2} we showed that $\pi(\phi,\nu|\y)$ can be expressed in terms of an arbitrary prior distribution.  To proof the propriety of the posterior distribution of $(\bbeta, \sigma^2,\phi, \nu)$ under the proposed prior, we need to ensure that
\begin{equation*}
\int_{0}^{+\infty} \int_{4 + \varepsilon}^{+\infty} A(\nu)L^{*}(\phi,\y)\pi(\nu,\phi)d \nu d\phi < \infty.
\end{equation*}

Under the some conditions, we have that,
\begin{enumerate}
\item if $\mathbf{1}$ is not a column of $\X$. Then,
$$ tr\left[\Phi\right]=O\left(\frac{d }{d \phi}\log\psi(\phi)\right);$$

\item if $\mathbf{1}$ is not a column of $\X$. Then,
$$ tr\left[\Phi\right]=O\left(\frac{\omega(\phi)}{\psi(\phi)}\frac{d }{d \phi}\log\left[\frac{\omega(\phi)}{\psi(\phi)}\right]\right).$$
\end{enumerate}

Next, we will proof the case when $\mathbf{1}$ is not a column of $\X$, the case when $\mathbf{1}$ is a column of $\X$ is analogous and will be omitted. Let $n_p=n-p$ and using the Stirling approximation for the trigamma function ($\Psi_1(x)=x^{-1}+ (2x^2)^{-1}$), we have that as $\nu \to \infty$ 
$$D \approx \frac{1}{4}\left \lbrace\frac{-2n_p^3\nu+2n_p^3-2n_p^2\nu^2-4n_p^2\nu+8n_p^2+8n_p\nu+8n_p}{\nu^2(n_p+\nu^2)(n_p+\nu+2)}\right\rbrace.$$

Note that the higher order for the numerator (in terms of $\nu$) is $2$, while the denominator is a polynomial of order $5$, so we have that $D=O(\nu^{-3})$ and, 
$$BCD= O(1)\left ( O(1)\mathbf{A}- O(1)tr^2\left[\Phi\right]\right ) O(\nu^{-3}).$$

Now, $A=O(1)\left(2tr\left[ \Phi^2\right]+tr^2\left[\Phi\right]\right)$. Then, as $\phi \to +\infty$, we have that
$$BCD=O(\nu^{-3})O\left(\left(\frac{d }{d \phi}\log\psi(\phi)\right)^2\right).$$

Analogously,
\begin{eqnarray*}
 B_{11}C_{11}B_{12}&=&\frac{2\nu n_p^2}{(\nu+n_p+2)^3(\nu-2)^2(\nu+n_p)^2} tr^2\left[\Phi\right]\\ &=&O(\nu^{-6})O\left(\frac{d }{d \phi}\log\psi(\phi)\right),
\end{eqnarray*}
and hence 
\begin{eqnarray} \nonumber
|I_1(\theta^{*})| &<&  BCD + B_{11}C_{11}B_{12}
=O(\nu^{-3})O\left(\left(\frac{d }{d \phi}\log\psi(\phi)\right)^2\right) + \\
&& O(\nu^{-6})O\left(\left(\frac{d }{d \phi}\log\psi(\phi)\right)^2\right)
\\ \nonumber
&=& O(\nu^{-3})O\left(\left(\frac{d }{d \phi}\log\psi(\phi)\right)^2\right).
\end{eqnarray}

Finally, as $\nu \to 4+\varepsilon$, all expressions involving $\nu$ are constant. So, we have that
$$ BCD + B_{11}C_{11}B_{12} = O\left(\left(\frac{d }{d \phi}\log\psi(\phi)\right)^2\right), \,\,\,{\it as}\,\,\, \phi \to +\infty$$ and the result follows.

\bibliographystyle{chicago}
\bibliography{biblio2}
\end{document}